\documentclass{agtart_a}
\pdfoutput=1

\usepackage{graphicx}

\title[Deforming hyperbolic 3--manifolds with boundary]
{On deformations of hyperbolic 3--manifolds with geodesic boundary}

\author{Roberto Frigerio}
\givenname{Roberto}
\surname{Frigerio}
\address{Dipartimento di Matematica\\
Universit\`a di Pisa\\\newline
Largo B~Pontecorvo 5\\
56127 Pisa\\
Italy}
\email{frigerio@mail.dm.unipi.it}
\urladdr{}

\volumenumber{6}
\issuenumber{}
\publicationyear{2006}
\papernumber{15}
\startpage{435}
\endpage{457}

\doi{}
\MR{}
\Zbl{}


\received{7 October 2005}
\revised{}
\accepted{20 February 2006}
\published{23 March 2006}
\publishedonline{23 March 2006}
\proposed{}
\seconded{}
\corresponding{}
\editor{}
\version{}

\arxivreference{math.GT/0504116}

\AtBeginDocument{\let\bar\wbar\let\tilde\wtilde\let\hat\what\let\overline\wbar}

\makeatletter
\def\cnewtheorem#1[#2]#3{\newtheorem{#1}{#3}[section]
\expandafter\let\csname c@#1\endcsname\c@lemma}

\newtheorem{lemma}{Lemma}[section]
\cnewtheorem{teo}[lemma]{Theorem}
\cnewtheorem{prop}[lemma]{Proposition}
\cnewtheorem{cor}[lemma]{Corollary} 
\cnewtheorem{conj}[lemma]{Conjecture}

\theoremstyle{definition}
\cnewtheorem{defn}[lemma]{Definition}
\cnewtheorem{fact}[lemma]{Fact}
\cnewtheorem{quest}[lemma]{Question}
\cnewtheorem{probl}[lemma]{Problem}
\cnewtheorem{example}[lemma]{Example}
\cnewtheorem{rem}[lemma]{Remark} 

\makeautorefname{teo}{Theorem}
\makeautorefname{defn}{Definition}

\makeatother  

\newcommand{\matN}{\ensuremath {\mathbb{N}}}
\newcommand{\matR}{\ensuremath {\mathbb{R}}}

\newcommand{\matZ}{\ensuremath {\mathbb{Z}}}
\newcommand{\matC}{\ensuremath {\mathbb{C}}}

\newcommand{\matH}{\ensuremath {\mathbb{H}}}

\newcommand{\calR}{\ensuremath {\mathcal{R}}}

\newcommand{\calI}{\ensuremath {\mathcal{I}}}

\newcommand{\calS}{\ensuremath {\mathcal{S}}}

\newcommand{\calC}{\ensuremath {\mathcal{C}}}

\newcommand{\calF}{\ensuremath {\mathcal{F}}}

\newcommand{\calT}{\ensuremath {\mathcal{T}}}

\newcommand{\calW}{\ensuremath {\mathcal{W}}}

\newcommand{\Nbar}{{\wbar N}}
\newcommand{\Hhyp}{\matH_{\,{\rm hyp}}}

\newcommand{\hp}{\mathcal{H}_{+}^3}
\newcommand{\hm}{\mathcal{H}_{-}^3}
\newcommand{\Minkos}{\mathbb{M}^{\, 3+1}}

\newcommand{\diag}{\SelectTips{cm}{} \xymatrix@1}

\newcommand{\lp}{L^3_+}

\def\diag{\SelectTips{cm}{} \xymatrix@1}
\newcommand{\tetra}{\Delta}

\begin{document}

\begin{asciiabstract}
Let M be a complete finite-volume hyperbolic 3-manifold with compact
non-empty geodesic boundary and k toric cusps, and let T be a
geometric partially truncated triangulation of M.  We show that the
variety of solutions of consistency equations for T is a smooth
manifold or real dimension 2k near the point representing the unique
complete structure on M.  As a consequence, the relation between
deformations of triangulations and deformations of representations is
completely understood, at least in a neighbourhood of the complete
structure.  This allows us to prove, for example, that small
deformations of the complete triangulation affect the compact
tetrahedra and the hyperbolic structure on the geodesic boundary only
at the second order.
\end{asciiabstract}

\begin{htmlabstract}
Let <i>M</i> be a complete finite-volume hyperbolic 3&ndash;manifold
with compact non-empty geodesic boundary and k toric cusps, and let
<i>T</i> be a geometric partially truncated triangulation of <i>M</i>.
We show that the variety of solutions of consistency equations for
<i>T</i> is a smooth manifold or real dimension <i>2k</i> near the point
representing the unique complete structure on <i>M</i>.  As a consequence,
the relation between deformations of triangulations and deformations of
representations is completely understood, at least in a neighbourhood
of the complete structure.  This allows us to prove, for example, that
small deformations of the complete triangulation affect the compact
tetrahedra and the hyperbolic structure on the geodesic boundary only
at the second order.
\end{htmlabstract}

\begin{abstract}
Let $M$ be a complete finite-volume hyperbolic 3--manifold with
compact non-empty geodesic boundary and $k$ toric cusps, and let
$\mathcal{T}$ be a geometric partially truncated triangulation of $M$.
We show that the variety of solutions of consistency equations for
$\mathcal{T}$ is a smooth manifold or real dimension $2k$ near the
point representing the unique complete structure on $M$.  As a
consequence, the relation between deformations of triangulations and
deformations of representations is completely understood, at least in
a neighbourhood of the complete structure.  This allows us to prove,
for example, that small deformations of the complete triangulation
affect the compact tetrahedra and the hyperbolic structure on the
geodesic boundary only at the second order.
\end{abstract}

\maketitle

The idea of constructing hyperbolic structures on manifolds by 
suitably gluing to each other geodesic polyhedra dates back to
Thurston~\cite{Thu:bibbia}. In the setting of cusped 
manifolds one employs ideal tetrahedra,
which are parameterized by complex numbers, and tries to solve
hyperbolicity equations.
In~\cite{FriPe} (written jointly with Petronio)
we explained how this
approach can be adapted to
the case of non-empty geodesic boundary:
in the bounded case 
one has to consider truncated tetrahedra, whose parameterization is 
more complicated, but basically the whole scheme extends.

The conditions under which  
a gluing of truncated tetrahedra
defines a non-singular hyperbolic metric are encoded 
by \emph{consistency equations}, 
while \emph{completeness equations}  
translate the conditions ensuring that such a 
metric is complete. For our purposes it is crucial
to control the number of consistency equations, and
this is the reason why the equations described here are quite different
from those introduced in~\cite{FriPe}.
The set of solutions of consistency equations naturally 
provides a deformation space for finite-volume
hyperbolic structures with geodesic boundary
on a fixed $3$--manifold.
Building on classical results in cohomology theory of
representations, 
we prove that the complete structure
is a smooth point of this deformation space and we explicitly 
describe local coordinates around it. This allows us to
give a proof of Thurston's hyperbolic Dehn filling Theorem
which applies to all
the hyperbolic manifolds with geodesic boundary which admit
a \emph{good} geometric triangulation (see \fullref{goodtria:def}). 
There is strong evidence that any complete
finite-volume hyperbolic $3$--manifold with geodesic boundary 
should admit
a good triangulation (see the discussion preceding \fullref{conj:conj});
moreover, any such manifold admits
a \emph{partially flat} triangulation, and following Petronio and Porti~\cite{PetPor}
we could probably adapt our proof 
of Thurston's hyperbolic Dehn filling Theorem 
to deal 
also with this kind of triangulations. This would give a complete
and self-contained proof of the filling Theorem via deformation theory
of geometric triangulations. 

In the last section we show that small deformations of the 
solution representing the complete structure 
affect the compact tetrahedra and the hyperbolic structure
on the geodesic boundary only at the second order.
These results are extensively 
used in~\cite{Fri:preprint} for studying small deformations 
in infinitely many concrete examples.

It is maybe worth mentioning that 
the deformation variety defined by the consistency
equations for a cusped $3$--manifold without boundary has already been
studied by several authors (see eg, Neumann--Zagier~\cite{NeuZag},
Petronio--Porti~\cite{PetPor}). 
In particular, Choi has recently proved in~\cite{Choi}
that in the cusped empty-boundary case
the deformation variety is a smooth complex manifold at any point representing
a non-degenerate (ie, neither partially flat nor partially negatively
oriented) geodesic ideal triangulation.

The author is partially supported by the INTAS project
``CalcoMet-GT'' 03-51-3663.
 
\section{Triangulations and hyperbolicity equations}\label{trunc:tria:section}

Let $N$ be a complete 
finite-volume orientable hyperbolic 3--manifold with \emph{compact}
non-empty
geodesic boundary
(from now on we will usually 
summarize all this information saying just that \emph{$N$ is hyperbolic}).
It is well-known that $N$ consists of a compact portion containing
$\partial N$ together with
several cusps of the form
$T\times[0,\infty)$, where $T$ is the 
torus, so  $N$ admits
a natural compactification $\Nbar$
obtained by adding some boundary tori.
Since the components of $\partial N$ are totally geodesic,
they inherit a hyperbolic metric, and have therefore
negative Euler characteristic.

\subsection{Partially truncated tetrahedra}
A \emph{partially truncated tetrahedron} is a pair $(\tetra,\calI)$, where 
$\tetra$ is a tetrahedron and $\calI$ is a set of vertices of $\Delta$,
that will be called \emph{ideal vertices}.
In the sequel we will always
refer to $\tetra$ itself as a partially truncated tetrahedron, tacitly implying
that $\calI$ is also fixed. The \emph{topological realization} $\tetra^{\!\ast}$ 
of  $\tetra$
is obtained by removing from $\tetra$ the ideal 
vertices
and small open
stars of the non-ideal vertices. 
We call \emph{lateral hexagon} and \emph{truncation 
triangle} 
the intersection of $\tetra^{\!\ast}$ respectively with a face of $\tetra$ 
and with the link in $\tetra$ of a non-ideal vertex. The edges of the 
truncation triangles, which also belong to the lateral hexagons, are called 
\emph{boundary edges}, and the other edges of $\tetra^{\!\ast}$ are called
\emph{internal edges}. If $\tetra$ has ideal vertices, 
a lateral
hexagon of $\tetra^{\!\ast}$ may not be 
a hexagon, because some of its closed edges may be missing.

\begin{figure}[ht!]
\begin{center}
\input{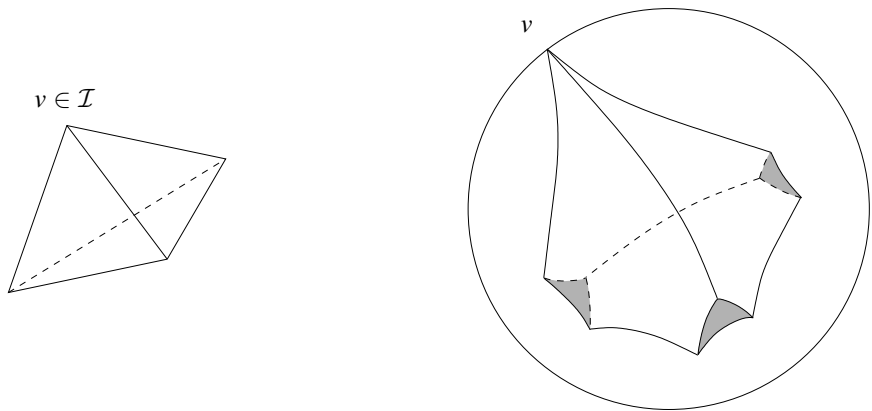}
\caption{A geometric tetrahedron
with one ideal vertex}\label{geomreal:fig}
\end{center}\end{figure}

A \emph{geometric realization} of $\tetra$ is an 
identification
of $\tetra^{\!\ast}$ with a convex polyhedron
in $\matH^3$ such that the truncation triangles
are geodesic triangles, the lateral hexagons are 
geodesic polygons with ideal vertices
corresponding to missing edges, and truncation 
triangles and lateral hexagons 
lie at right angles to each other.
An example of a 
geometric realization is shown in \fullref{geomreal:fig}, where
truncation triangles are shadowed.

\subsection{Triangulations}
Let $\Nbar$ be a compact orientable manifold and let
$N$ be obtained from $\Nbar$ by removing 
the toric components of $\partial\Nbar$.
We define a \emph{partially truncated triangulation} of $N$ to be 
a realization of $N$ as a gluing of some $\Delta\!^*$'s along
a pairing of the lateral hexagons induced by a simplicial
pairing of the faces of the $\Delta$'s.
When $N$ is endowed with a hyperbolic structure,
a partially truncated triangulation of $N$ is called \emph{geometric} 
if, for each tetrahedron $\Delta$ of the triangulation, the pull-back
to $\Delta^{\!\ast}$ of the Riemannian metric of $N$
defines a geometric realization of $\Delta$. Equivalently, 
the hyperbolic structure of $N$ 
should be obtained by gluing geometric realizations of the $\Delta$'s
along isometries of their lateral hexagons.

\begin{defn}\label{goodtria:def}
A partially truncated triangulation $\calT$ of an orientable
$3$--manifold $N$ is \emph{good} if any tetrahedron in $\calT$
has at most one ideal vertex.
\end{defn}

Kojima proved in~\cite{Koj:proc} that every hyperbolic $N$ has a
\emph{canonical} decomposition into partially truncated
\emph{polyhedra}, rather than tetrahedra.  Any polyhedron in a
canonical decomposition has at most one ideal vertex, so
triangulations arising as subdivisions of Kojima decompositions are
examples of good triangulations.  In the vast majority of cases the
Kojima decomposition actually consists of tetrahedra, or at least can
be subdivided into a geometric partially truncated triangulation.  For
instance, it is proved in Frigerio--Martelli--Petronio~\cite{FriMaPe2}
that there exist exactly $5192$ hyperbolic manifolds with non-empty
geodesic boundary which can be (topologically) triangulated by at most
four partially truncated tetrahedra: their Kojima decomposition can
always be subdivided into a triangulation, and is itself a
triangulation in $5108$ cases.  These facts strongly support the
following:

\begin{conj}\label{conj:conj}
Any hyperbolic $N$ with non-empty geodesic boundary 
admits a good geometric triangulation.
\end{conj}

\subsection{Moduli for partially truncated tetrahedra}
The following result implies that the dihedral angles can
be used as moduli for geometric tetrahedra. 
\begin{teo}\label{moduli:teo}
Let $\Delta$ be a partially truncated tetrahedron and let 
$\Delta\!^{(1)}$ be the
set of edges of $\Delta$. The geometric realizations of $\Delta$ are
parameterized up to isometry 
by the dihedral angle assignments $\theta\co \Delta\!^{(1)}\to(0,\pi)$ such that
for each vertex $v$ of $\Delta$, 
if $e_1,e_2,e_3$ are the edges that emanate from $v$, then
$\theta(e_1)+\theta(e_2)+\theta(e_3)$ is equal to $\pi$ for ideal $v$
and less than $\pi$ for non-ideal $v$.
\end{teo}

\begin{figure}[ht!]
\begin{center}
\input{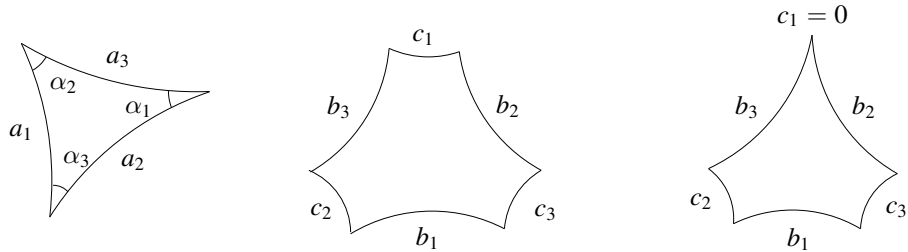}
\caption{A triangle, a right-angled hexagon and a pentagon
with four right angles and an ideal vertex}\label{trigo:fig}
\end{center}
\end{figure}

Having introduced moduli for geometric tetrahedra, our next task is 
to determine, given a triangulated manifold, which values of moduli 
define a global hyperbolic structure on the manifold.
The following well-known hyperbolic trigonometry formulae will prove useful
later:

\begin{lemma}\label{trigo:lemma}
With notation as in \fullref{trigo:fig} we have
\begin{eqnarray*}
\cosh a_1={(\cos \alpha_2\cdot\cos \alpha_3+\cos\alpha_1)}/{(\sin
\alpha_2\cdot\sin\alpha_3)},\label{cos:rule}\\
\cosh b_1={(\cosh c_2\cdot\cosh c_3+\cosh c_1)}/{(\sinh
c_2\cdot\sinh c_3)}\label{hexa:rule}.
\end{eqnarray*}
\end{lemma}

\begin{figure}[ht!]
\begin{center}
\input{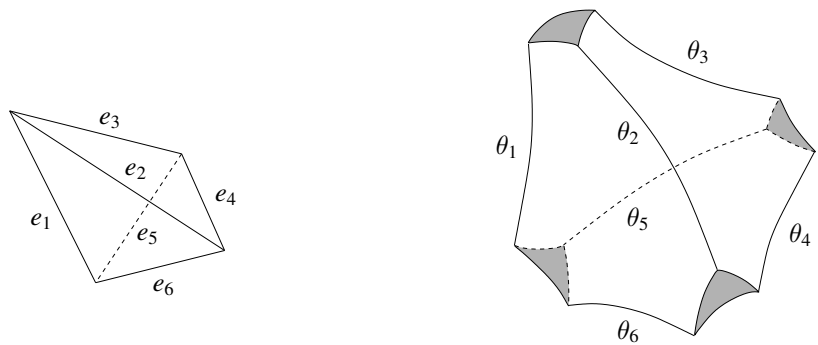}
\caption{Notation for edges and dihedral angles 
of a truncated tetrahedron}\label{tetra:notation:fig}
\end{center}\end{figure}

Let now $\Delta$ be a partially truncated tetrahedron 
with edges $e_1,\ldots,e_6$
as in \fullref{tetra:notation:fig}. 
We fix a geometric realization $\theta$ of $\Delta$
determined by the dihedral angles $\theta_i=\theta(e_i)$
for $i=1,\ldots,6$ and we denote by $L^{\theta}$ the
length with respect to this realization.
The boundary edges of the lateral hexagons of $\Delta$
correspond to the pairs of distinct non-opposite 
edges $\{e_i,e_j\}$, and will be denoted by $e_{ij}$.
Note that $e_{ij}$ disappears towards infinity, so it has length $0$, when the 
common vertex of $e_i$ and $e_j$ is ideal. 
\fullref{trigo:lemma}
readily implies
\begin{equation}\label{boundary:edge:length:formula}
\cosh L^\theta(e_{12})=({\cos \theta_1 \cdot \cos \theta_2 + \cos \theta_3})/
({\sin \theta_1 \cdot \sin \theta_2}).
\end{equation}
Note that this result is correct also when
the common end of $e_1$ and $e_2$ is ideal.

Turning to the length of an internal edge, we note that the edge 
is an infinite half-line or an infinite line when one or both 
its ends are ideal, respectively.
Otherwise the length is computed 
using \fullref{trigo:lemma}.
With notation as in \fullref{tetra:notation:fig}, and
defining $v_{ijk}$ as the vertex from which the edges
$e_i,e_j,e_k$ emanate, we set:
\begin{eqnarray*}
\begin{array}{rcl}
c^\theta(e_1)\!\!&=&\!\!\cos \theta_1\cdot \left( \cos \theta_3\cdot \cos\theta_6
        +\cos \theta_2\cdot \cos\theta_5\right)\\
        & & + \cos \theta_2\cdot \cos\theta_6
        +\cos \theta_3\cdot \cos\theta_5 + \cos\theta_4\cdot \sin^2 \theta_1;
\end{array}\\
d^\theta(v_{123})\,=\,2\cos \theta_1\cdot \cos\theta_2\cdot \cos\theta_3+ \cos^2 \theta_1
+\cos^2 \theta_2+\cos^2 \theta_3 -1.
\end{eqnarray*}

\begin{prop}\label{internal:edge:length:prop}
$d^\theta(v_{123})=0$ if and only if
the vertex $v_{123}$ is ideal.
If $v_{123}$ and $v_{156}$ are both non-ideal then
\begin{equation}\label{internal:edge:length:formula}
\cosh L^\theta(e_1)={c^\theta(e_1)}\,\Big/\,{\sqrt{d^\theta(v_{123})
\cdot d^\theta(v_{156})}}.
\end{equation}
\end{prop}

\begin{rem}\label{lengths:rem}
Let $\Delta$ be a partially truncated tetrahedron without
ideal vertices. Then the geometric realizations of $\Delta$ are
parameterized by the 
lengths of the internal edges. In fact, the map that associates 
to the dihedral
angles of a geometric realization of $\Delta$ the lengths of its 
internal edges is a diffeomorphism between open subsets
of $\matR^6$.
\end{rem}

\subsection{Conditions for geometric gluing}
Let $N$ be obtained from an orientable compact $\Nbar$ by removing
all the tori in $\partial\Nbar$ and let $\calT$ be a partially truncated
triangulation of $N$. Let also $\theta$ be a geometric realization
of the tetrahedra in $\calT$ and denote by $L^\theta$
the length with respect to this realization. 
We now describe the conditions under which the realization
$\theta$ defines a hyperbolic structure on the whole of $N$.
For our purposes it will be sufficient
to deal only
with \emph{good\/} triangulations, so we assume from now on
that $\calT$ is good. The general case is treated 
in~\cite{FriPe}.

In order to define a global hyperbolic structure on $N$,
the tetrahedra of $\calT$ must satisfy two obvious necessary 
conditions, which in fact are also sufficient. Namely, 
we should be able
to glue the lateral hexagons by isometries, and
we should have a total dihedral angle of $2\pi$ around each 
edge of the manifold.
The first condition ensures that the hyperbolic structure defined
by $\theta$ on the complement of the $2$--skeleton of $\calT$
extends to the complement of the $1$--skeleton.  Since $\calT$ is good,
the second one ensures that the structure glues up without singularities 
also along the edges. 
The second condition is directly expressed in terms of moduli, and we
will explain in a moment how to translate the first one 
into an equation on dihedral angles.

\begin{rem}\label{nogood:rem}
If $\calT$ were not good, requiring the dihedral angles
around each edge to sum up to $2\pi$ would not be sufficient
to obtain a non-singular hyperbolic metric on the $1$--skeleton of $\calT$.
The point is that when some geometric tetrahedra are arranged one after
the other around an edge $e$ with two ideal endpoints,  
the first face of the first tetrahedron and the second face of the last 
tetrahedron may overlap without coinciding. Namely, the isometry which pairs
these two faces may be a translation along $e$ instead of being the identity.
Of course the isometry has to be the identity if at least one endpoint
of $e$ is not ideal.
\end{rem}

\subsection{Exceptional hexagons}
It is easily seen  that
gluings match ideal vertices to each other, because 
these notions are part of the 
initial topological information about a triangulation.
When a pairing glues two \emph{compact} lateral hexagons, 
to be sure that the gluing is an isometric one
we may equivalently require the lengths of the internal edges or
those of the boundary edges to match under the gluing.
On the other hand, since $\calT$ is a good triangulation, a non-compact 
lateral hexagon $F$
is actually a pentagon with four right angles and an ideal vertex:
we shall say in this case that $F$ is an \emph{exceptional} lateral hexagon.
By \fullref{trigo:lemma}, 
the isometry class of an exceptional lateral hexagon
is determined by the lengths of its boundary edges. However,
in order to end up with a non-redundant set of consistency equations,
it is convenient to 
find an alternative approach to moduli 
for  exceptional hexagons.
\begin{figure}[t]
\begin{center}
\input{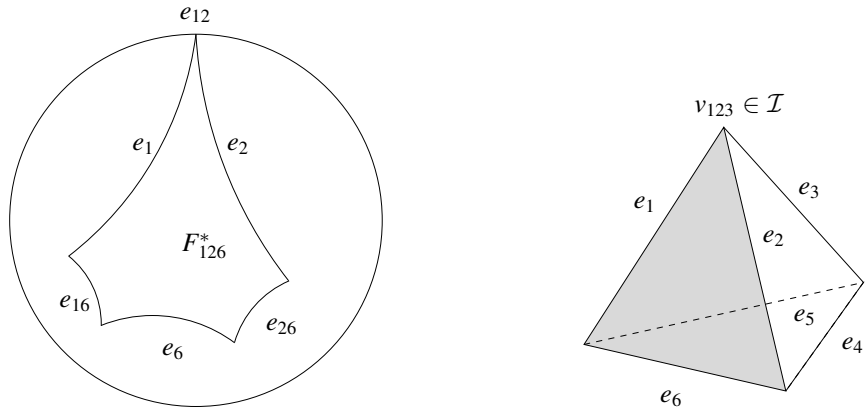}
\caption{An exceptional hexagon}\label{special:hexagon:fig}
\end{center}\end{figure}
To this aim we need now to be slightly more
careful about orientation than we have been so far. 
Namely, we choose on the tetrahedra an orientation compatible
with a global orientation of the manifold. As a result
also the lateral hexagons
have a fixed orientation, and the gluing maps
reverse the orientation of the hexagons.

So, let us consider an exceptional hexagon $F^*_{126}$ 
as in \fullref{special:hexagon:fig}, and  recall that the hexagon is
oriented and embedded in $\matH^3$ by $\theta$.
We consider the horospheres $O_1$ and
$O_2$ centred at $e_{12}$ and passing through the 
non-ideal ends of $e_1$ and $e_2$ respectively.
We define $\sigma^\theta(F_{126})$ to be
$\pm{\rm dist}(O_1,O_2)$, the sign being positive
if $e_2,e_{12},e_1$
are arranged positively on $\partial F^*_{126}$
and $O_1$ is contained in the horoball bounded
by $O_2$, or if $e_2,e_{12},e_1$ are arranged negatively
on $\partial F^*_{126}$ and $O_2$ is contained in the horoball
bounded by $O_1$, and 
negative otherwise. Together with 
equation~\ref{boundary:edge:length:formula}, the following proposition
allows to compute $\sigma$ in terms of the dihedral angles.

\begin{prop}\label{compute:sigma:prop}
We have 
\begin{equation}\label{compute:sigma:formula}
\sigma^\theta(F_{126})=\ln (\sinh L^\theta (e_{16})/\sinh L^\theta (e_{26})).
\end{equation}
\end{prop}
\begin{proof}
Let $\alpha_1, \alpha_2$ and $R, R_1, R_2$ be 
the angles and lengths shown in \fullref{sigma:comp:fig}.
An easy computation in the upper half-plane model of
$\matH^2$ shows that $L^\theta (e_{i6})=
\ln \cot (\alpha_i/2)$ for $i=1,2$.
Moreover we have $R=R_1 \cdot \tan \alpha_1= R_2 \cdot\tan \alpha_2$, so
$\exp \sigma^\theta(F_{126})=R_1/R_2=\tan \alpha_2/\tan \alpha_1$.
\begin{figure}[t]
\begin{center}
\input{\figdir/hexa.pstex_t}
\caption{Notation for the proof of 
\fullref{compute:sigma:prop}}\label{sigma:comp:fig}
\end{center}\end{figure}
Now for $i=1,2$ we have $\tan \alpha_i=
2\cot (\alpha_i /2)/(\cot^2 (\alpha_i /2) -1)= 
2\exp L^\theta (e_{i6})/(\exp (2 L^\theta (e_{i6})) -1)$.
Combining these equations we finally get
$\exp \sigma^\theta(F_{126})=
\sinh L^\theta (e_{16})/\sinh L^\theta (e_{26})$. 
\end{proof}

We now define $\ell^\theta(F_{126})$ to be
the length of $e_6$.
The next proposition shows that the functions 
$\sigma$ and $\ell$ provide a parameterization
of isometry classes of exceptional hexagons.

\begin{prop}\label{sigma:works:prop}
Let $F$ and $F'$ be paired exceptional lateral hexagons.
Their pairing can be realized by an isometry if and only
if $\sigma^\theta(F)+\sigma^\theta(F')=0$
and $\ell^\theta(F)=\ell^\theta(F')$.
\end{prop}
\begin{proof}
We concentrate on the ``if'' part of the statement, the ``only if''
part being obvious.
Let $f_1,f_2$ and $f'_1, f'_2$ be the boundary edges of
of $F$ and $F'$ respectively, and assume that the orientation-reversing
pairing between $F$ and $F'$ glues $f_1$ to $f'_2$
and $f'_1$ to $f_2$. Since  
$\ell^\theta(F)=\ell^\theta(F')$,
\fullref{trigo:lemma} easily implies that 
$L(f_1)\leqslant L(f'_2)$ if and only if
$L(f_2)\geqslant L(f'_1)$. Moreover the assumption
that $\sigma^\theta(F)=-\sigma^\theta(F')$
and \fullref{compute:sigma:prop} give
$L(f_1)\leqslant L(f'_2)$ if and only if
$L(f_2)\leqslant L(f'_1)$. This forces
$L(f_1)= L(f'_2)$ and
$L(f_2)=L(f'_1)$, whence the conclusion.
\end{proof}

\begin{rem}\label{sigma:rem}
Let $\Delta$ be a partially truncated tetrahedron and suppose
that $v$ is the unique ideal vertex of $\Delta$. Let 
$F_1,F_2,F_3$ be the faces of $\Delta$ incident to $v$ and 
for $i=1,2,3$ let $e_i$ be the edge of $F_i$ not containing
$v$. The isometry classes of the geometric realizations
of $\Delta$ are parameterized by the lengths of the 
$e_i$'s and the values taken by $\sigma$ on the $F_i$'s.
More precisely, the map that associates to any geometric
realization $\theta$ of $\Delta$ the point
$(L^\theta (e_1),L^\theta (e_2),L^\theta (e_3),\sigma^\theta (F_1),
\sigma^\theta (F_2),\sigma^\theta (F_3))$ defines a diffeomorphism
between open subsets of two affine hyperplanes of $\matR^6$.
\end{rem} 

\subsection{Consistency equations}
Recall now that we are considering a 
candidate hyperbolic 
$3$--manifold
$N$ endowed with a good triangulation $\calT$,
and that we have fixed
a geometric realization $\theta$ 
of the tetrahedra in $\calT$. 
The above discussion implies the following:

\begin{teo}\label{consistency:teo:second} 
The parameterization
$\theta$ defines on $N$ a hyperbolic structure
with  geodesic boundary if and only if the following 
conditions hold:
\begin{enumerate}
\item the total dihedral angle along any edge of $\calT$ in $N$
is equal to $2\pi$;
\item $L^\theta(e)=L^\theta(e')$ for
all pairs $(e,e')$ of matching compact internal edges;
\item
$\sigma^\theta(F)+\sigma^\theta(F')=0$ for all pairs 
$(F,F')$ of matching exceptional hexagons.
\end{enumerate}
\end{teo}

By Propositions~\ref{internal:edge:length:prop},~\ref{compute:sigma:prop}, 
conditions~(1), (2), (3) 
of \fullref{consistency:teo:second} translate into a set
$\calC(\calT)$ of smooth equations on $\theta$,
which are called 
\emph{consistency equations}.
Our next task is to compare the number of
equations in $\calC(\calT)$ 
with the dimension of the moduli space of
geometric realizations of the tetrahedra of $\calT$.
  
Let $c$ (resp.~$p$) be the number of compact (resp.~non-compact)
tetrahedra of $\calT$. By \fullref{moduli:teo},
the number of parameters for the geometric realizations  
of $\calT$ is equal to $6 c+ 5 p=6t-p$, where $t$ is the total
number of tetrahedra of $\calT$.
Let $e_1,\ldots, e_l$ be the edges of $\calT$ without ideal endpoints,
and for $i=1,\ldots, l$ let $x_i$ be the valence of $e_i$, ie, the
number of tetrahedra of $\calT$ incident to $e_i$, with multiplicity.
Of course we have $\sum_{i=1}^l x_i=6c+ 3p$.
If we denote by $h$ the number of 
edges of $\calT$ with exactly one ideal endpoint, then
the union of the boundary tori of $\Nbar$
admits a triangulation with $h$ vertices and $p$ triangles,
so $h=p/2$.
Let us now turn to the number of equations. 
\fullref{consistency:teo:second} determines $l+h$ equations
arising from conditions~(1),
$\sum_{i=1}^l (x_i-1)$ equations
arising from conditions~(2) and $3p/2$ equations
arising from conditions~(3). Since 
\begin{equation*}
\begin{array}{lllll}
& & l+h+\left(\sum_{i=1}^l (x_i-1)\right)
+3p/2 &=& h+\left(\sum_{i=1}^l x_i\right) +3p/2\\ &=&
p/2+6c+3p+3p/2 &=& 6c+5p,
\end{array}
\end{equation*} 
we can conclude that
the number of equations in $\calC(\calT)$ is equal to the 
dimension of the moduli space of
geometric realizations of the tetrahedra of $\calT$.

\subsection{Reducing the number of equations}
When $N$ has cusps, some equations in $\calC(\calT)$
turn out to be redundant. Let $T$ be a fixed
toric component of $\partial \Nbar$ and let $j$ be the number
of tetrahedra asymptotic to $T$.
Since any such tetrahedron contributes to the triangulation
of $T$ with a Euclidean triangle,
the sum of all the dihedral angles along all the edges
of $\calT$ incident to $T$ is equal to $j\pi$. 
So if we require condition~(1) of \fullref{consistency:teo:second}
to hold for \emph{all but one} edge incident to $T$,
then the same condition is automatically satisfied also along the remaining
edge. This allows us to discard from $\calC(\calT)$
one equation for each cusp of $N$. 

Moreover, let $\Delta$ be a partially truncated tetrahedron
with an ideal vertex $v$ incident to $T$, and let $F_1,F_2,F_3$ be the 
faces of $\Delta$ incident to $v$. By the very definition of $\sigma$
it follows that $\sigma(F_1)+\sigma(F_2)+\sigma(F_3)=0$. 
This implies that if $\calF$ is the set of all the faces incident to $T$
of tetrahedra of $\calT$, then we have 
$\sum_{F \in\calF} \sigma (F)=0$. So if we require condition~(3) of 
\fullref{consistency:teo:second} to hold for \emph{all but one} pair 
of matching exceptional hexagons incident to $T$, then the same condition
also holds for the remaining pair. This means that another equation of $\calC(\calT)$
for each cusp of $N$ can be discarded.
Suppose that $N$ compactifies to an orientable $\Nbar$ 
with $k$ boundary tori.
The above discussion is summarized by the following:

\begin{prop}\label{reduce:eq:prop}
We can discard $2k$ equations
from $\calC(\calT)$ thus obtaining an equivalent set
$\calC^\ast(\calT)$ of $n-2k$ equations, where $n$
is the dimension of the moduli space of geometric
realizations of the tetrahedra of $\calT$.
\end{prop}

We have seen in the preceding subsection that
the moduli space of geometric realizations
of the tetrahedra of $\calT$  
is given by a subset
$\calW$ of $\matR^{6t}$, where $t$ is the number of tetrahedra
of $\calT$. More precisely, $\calW$ is an open convex subset
of an affine subspace of dimension $6t-p$, where $p$
is the number of non-compact tetrahedra in $\calT$.
Recall that $l$ is the number of compact edges (considered
as subsets of $N$) of $\calT$. The above computation implies that
equations in $\calC^\ast (\calT)$ corresponding to conditions~(1)
of \fullref{consistency:teo:second} take the form
$A (x)=0$, where $A\co \calW\to\matR^{l+p/2-k}$ is an affine map,
while equations corresponding to conditions~(2)
and~(3)
take the form
$F (x)=0$, where $F\co \calW\to\matR^{6t-l-(3/2)p-k}$ is constructed
from
formulae~\ref{internal:edge:length:formula},~\ref{compute:sigma:formula}, 
and is therefore smooth.
From now on we denote by $\Omega(\calT)=F^{-1} (0)\cap A^{-1} (0)
\subset\calW\subset
\matR^{6t}$ the set of
solutions of consistency equations $\calC^\ast(\calT)$.

\subsection{Completeness}\label{completeness:subsection}
Let
$T_1,\ldots,T_k$ be the boundary tori of $\Nbar$. From now on 
we denote by $\mu_i,\lambda_i$ a fixed basis of
$H_1 (T_i;\matZ)\cong\pi_1 (T_i)$, $i=1,\ldots,k$. 
Any point in $\Omega(\calT)$
naturally defines an ${\rm Aff} (\matC)$--structure
on $T_i$ (see eg, \cite{BenPet:book,Fri:tesi}).
For $x\in\Omega (\calT)$, we denote by 
$a_i (x)\in\matC$ (resp.~by $b_i (x)\in\matC$) the linear component
of the holonomy of $\mu_i$ (resp.~of $\lambda_i$) corresponding
to the ${\rm Aff} (\matC)$--structure defined by $x$ on $T_i$.
It is well-known that the hyperbolic
structure defined by $x$ on $N$ induces a complete metric on
the $i^{\rm th}$ cusp of $N$ if and only if $a_i (x)=b_i (x)=1$.
Moreover, one can explicitly compute $a_i$ and $b_i$
in terms of the dihedral angles as follows.

Let $\Delta$ be a tetrahedron in $\calT$, let $v$ be an ideal
vertex of $\Delta$ and $L_x(v)$ be the (similarity class of the)
Euclidean triangle obtained by intersecting
the geometric realization of $\Delta$ parameterized by $x$
with a small horosphere centred at $v$.
The tetrahedron being
oriented, this triangle is also oriented, so, once a vertex $p$ 
of $L_x(v)$ is fixed, we can associate to the similarity structure of $L_x(v)$
the unique complex number $z_x(L(v),p)$ 
such that $L_x(v)$ is carried to the Euclidean 
triangle with vertices $0,1,z_x(L(v),p)$
by an orientation-preserving similarity sending $p$ to $0$. 
Suppose that $e_1,e_2,e_3$ are the internal edges emanating
form $v$, and that they are positively arranged around $v$.
If $p=L(v_{123})\cap e_1$, then  
\begin{equation*}
z_x(L(v_{123}),p)=
({\sin \theta_2}/{\sin \theta_3})\cdot e^{i\theta_1}.
\end{equation*}

If $\gamma$ is an oriented simplicial loop
on $T_i$ and $q$ is a vertex of $\gamma$, then the set
of all triangles touching $\gamma$ in $q$ and lying on the right 
of $\gamma$ is well-defined and will be denoted by $R(\gamma,q)$. Moreover, we
shall denote by $V(\gamma)$ the set of vertices of $\gamma$.
Let $\hat{\mu}_i,\hat{\lambda}_i$ be simplicial loops
on $T_i$ representing $\mu_i,\lambda_i$.
The following result is proved in~\cite{Thu:bibbia,BenPet:book}.

\begin{prop}\label{complete:form:prop}
We have
\begin{equation*}
\begin{array}{lll}
a_i (x)&=&(-1)^{\# V(\hat{\mu}_i)} \cdot \prod_{q\in V(\hat{\mu}_i)}
\prod_{T\in R(\hat{\mu}_i,q)} z_x(T,q),\\
b_i (x)&=&(-1)^{\# V(\hat{\lambda}_i)} \cdot \prod_{q\in V(\hat{\lambda}_i)}
\prod_{T\in R(\hat{\lambda}_i,q)} z_x(T,q).
\end{array}
\end{equation*}
\end{prop}

As a consequence of Mostow--Prasad's rigidity Theorem
for hyperbolic manifolds with geodesic boundary~\cite{FriPe,Fri}
we get the following:

\begin{teo}\label{uniqueness:teo}
There exists at most one point in $\Omega (\calT)$
that defines on $N$ a complete hyperbolic structure
with geodesic boundary.
\end{teo}
\begin{proof}
See~\cite{FriPe}.
\end{proof}

\section{Smoothness at the complete structure}\label{smoothness:section}
Suppose now that $x_0$ is the unique point in $\Omega(\calT)$
which defines on $N$ a complete hyperbolic structure. For $i=1,\ldots,k$ and 
$x\in\Omega(\calT)$ let us define
\begin{equation*}
u_i (x)=\ln a_i (x),\qquad
v_i (x)=\ln b_i (x),
\end{equation*}
where 
$\ln$ is the branch of the complex logarithm defined on
$\{z\in\matC:\ \Re (z)>0\}$ such that $\ln 1=0$.
Since a non-trivial parabolic isometry does not commute with
a non-trivial orientation-preserving non-parabolic isometry,
we have the following:

\begin{prop}\label{ubasta:prop}
In a neighbourhood of $x_0$
in $\Omega (\calT)$ we have 
$u_i (x)=0$ $\Leftrightarrow$  $v_i (x)=0$
$\Leftrightarrow$ the  
hyperbolic structure defined by $x$ on the $i^{\rm th}$ cusp of $N$
is complete.
 \end{prop}

Let
$F\co \calW\to\matR^{6t-l-(3/2)p-k}$, 
$A\co \calW\to\matR^{l+p/2-k}$
be the smooth functions previously defined such that
$\Omega (\calT)=F^{-1}(0)\cap A^{-1}(0)$. 
We now set 
\[
G\co \calW \to \matR^{6t-p-2k}\times \matC^k,\quad
G(x)= (F(x), A(x), u_1(x),\ldots,u_k(x)).
\]
By \fullref{uniqueness:teo}
and \fullref{ubasta:prop} 
we have $G^{-1} (0)=\{x_0\}$.
This section 
is entirely devoted to the proof of our main result:

\begin{teo}\label{smoothness:teo}
We have ${\rm Ker}\ dG_{x_0}=\{0\}$. Thus:
\begin{enumerate}
\item 
$G$ induces a diffeomorphism of an open
neighbourhood of $x_0$ in $\calW$ onto an open neighbourhood 
of $0$ in $\matR^{6t-p-2k}\times \matC^k$;
\item
$\Omega(\calT)$ is a smooth manifold
of real dimension $2k$ near $x_0$;
\item
the map
\[
u\co \Omega(\calT)\to \matC^k,\quad u(x)=(u_1(x),\ldots,u_k(x))
\]
induces a diffeomorphism  of an open  neighbourhood of
$x_0$ in $\Omega (\calT)$
onto an open neighbourhood of $0$ in $\matC^k$.
\end{enumerate}
\end{teo}

Any point in $\Omega_{\rm c} (\calT):=F^{-1}(0)$ defines
a cone structure on $N$ having cone singularities along
the edges of $\calT$.
The following immediate consequence of \fullref{smoothness:teo}
was also proved  in~\cite{HodKer,Luo}.

\begin{cor}\label{smoothness3:cor}
If $N$ is compact
then $\Omega_{\rm c} (\calT)$ is parameterized in a neighbourhood
of $x_0$ by
the cone angles along the edges of $\calT$.
\end{cor}

\subsection{Deforming cone structures}
We begin with the following:

\begin{prop}\label{smooth:cone:prop}
The tangent map $dF_{x_0}\co  T_{x_0} (\calW)\to 
\matR^{6t-l-(3/2)p-k}$ has maximal rank, so
$\Omega_{\rm c} (\calT)$ is a manifold of real dimension
$l+p/2+k$ near $x_0$.
\end{prop}
\begin{proof}
By Remarks~\ref{lengths:rem},~\ref{sigma:rem},
the lengths of the 
compact internal edges and 
the values taken by $\sigma$ on the exceptional
lateral hexagons of the tetrahedra of $\calT$ 
provide smooth coordinates on $\calW$.
It is easily seen that
with respect to these coordinates the map
$F$ is affine and has maximal rank at $x_0$
(whence at any point of $\calW$).
\end{proof}

Let now $\alpha\co (-\varepsilon,\varepsilon)\to\Omega_{\rm c}
(\calT)$ be a smooth arc with ${\alpha}(0)=x_0$.
We will study the deformation associated to $\alpha$ using 
tools from the 
cohomology theory of representations: notation is as in 
\fullref{cohomology:app}, where we give some basic definitions
and results.

From now on we denote by $N'$ the non-compact manifold obtained
by drilling from $N$ all the edges of $\calT$. 
For any $t\in
(-\varepsilon,\varepsilon)$ the point $\alpha(t)$ determines a 
smooth hyperbolic structure $M_t$ on $N'$, whose completion
gives a hyperbolic cone structure on $N$. 
First of all we describe how to deduce the shape of the geometric
tetrahedra corresponding to the point $\alpha(t)$ just
from the geometric structure $M_t$ on $N'$. 
To this end we fix for $t\in(-\varepsilon,\varepsilon)$
a developing map $D_t\co \widetilde{N}'\to\matH^3$ with 
associated holonomy representation $\rho_t\co \pi_1 (N')\to
{\rm PSL}(2,\matC)$ (note that we can choose $D_t$ and $\rho_t$ to
vary smoothly with $t$).
Recall that $N$
compactifies
to a manifold $\Nbar$ with $k$ boundary tori, and denote by
$\overline{N}'$ the non-compact manifold obtained by drilling
from $\Nbar$ the \emph{closed} properly embedded arcs corresponding
to the edges of $\calT$. For $i=1,\ldots,k$ let $T'_i$ be the 
punctured torus $T'_i=T_i\cap \Nbar'$, and denote by $\calS_{\rm par}$
the family of all the boundary components of the universal
covering of $\overline{N}'$ projecting to some $T'_i$. 
We say that $P\subset\pi_1(N')\cong\pi_1(\overline{N}')$
is a \emph{parabolic peripheral subgroup} of $\pi_1(N')$ if
$P$ is the stabilizer of some boundary component $S\in\calS_{\rm par}$ 
of the universal covering of
$\overline{N}'$ (so
$\rho_0 (P)$ is a $\matZ+\matZ$
parabolic subgroup of ${\rm PSL}(2,\matC)$). 
Let $\partial_{\rm gd} \overline{N}'$
be the portion of $\partial \overline{N}'$ 
corresponding to the geodesic boundary of $N$, ie, let
$\partial_{\rm gd} \overline{N}'=\partial \overline{N}'\cap
\partial N$, and denote by $\calS_{\rm gd}$
the family of all the boundary components of the universal covering of 
$\overline{N}'$ projecting to some component of
$\partial_{\rm gd} \overline{N}'$.   
We say that $P\subset\pi_1(N')$
is a \emph{Fuchsian peripheral subgroup} of $\pi_1(N')$ if
$P$ is
the stabilizer of some component $S\in\calS_{\rm gd}$ (so
$\rho_0 (P)$ is a Fuchsian subgroup
of ${\rm PSL}(2,\matC)$).

\subsection{Dual vectors to planes and horospheres}
If $S\in \calS_{\rm gd}$, by construction the image of $S$ under 
the developing map $D_t$ is contained in a totally geodesic
immersed surface in $\matH^3$. It is easily seen that such a surface 
must in turn be contained in a geodesic plane $\hat{S} (t)$ 
of $\matH^3$. We now need to associate to each $\hat{S} (t)$
a suitable \emph{ultra-ideal} point $b_S (t)$, which will be called
the \emph{dual} point of $\hat{S} (t)$. Such point naturally lies in  
4--dimensional Minkowsky space,
so we fix some notation about this space.

We denote by $\Minkos$ the space $\mathbb{R}^4$ 
with coordinates $x_0,x_1,x_2,x_3$ 
endowed with the Lorentzian
inner product ${\langle x,y\rangle=-x_0 y_0+x_1 y_1+x_2 y_2+x_3 y_3}$. 
We set 
\begin{eqnarray*}
\hm&=&\{x\in\Minkos:\ \langle x,x \rangle=-1,\ x_0>0\},\\
\hp&=&\{x\in\Minkos:\ \langle x,x \rangle=1\},\\
\lp&=&\{x\in\Minkos:\ \langle x,x \rangle=0,\ x_0>0\}.
\end{eqnarray*}
We recall that $\hm$ is the upper sheet of the two-sheeted hyperboloid, 
and that
$\langle\,\cdot\,,\,\cdot\,\rangle$ restricts to a Riemannian metric
on $\hm$. With this metric,
$\hm$ is the so-called \emph{hyperboloid model} $\Hhyp^3$ of hyperbolic space.
The one-sheeted hyperboloid $\hp$ turns out to have a bijective
correspondence with the set of hyperbolic half-spaces in $\Hhyp^3$. Given
$b\in\hp$, the corresponding half-space, 
called the \emph{dual} of $b$, is given by
$\{v\in\hm:\ \langle v,b\rangle \leqslant 0\}$.
Similarly, the cone $\lp$ of future-oriented light-like vectors of $\Minkos$
corresponds to the set of horospheres in $\Hhyp^3$. The horosphere dual to
$u\in\lp$ is given by 
$\{v\in\hm:\ \langle v,u\rangle=-1\}$. 

Note now that for any $S\in\calS_{\rm gd}$ the set $D_t(\widetilde{N})$ 
locally lies on a definite side of $\hat{S} (t)$, so we can define
$b_S(t)$ to be the dual vector to the half-space that locally
contains $D_t(\widetilde{N})$ and is bounded by
$\hat{S} (t)$.   

When $S$ belongs to $\calS_{\rm par}$ a vector $b_S(t)\in\lp$ can also be defined
as follows: take an oriented edge $\widetilde{f}$ of $\widetilde{\calT}$ 
ending in $S$ and set $b_S(t)$ to be the unique point in $\lp$
with $x_0 (b_S(t))=1$ which projects to the endpoint
of $D_t (\widetilde{f})$ in $\partial \matH^3$. Such an endpoint exists
because $D_t (\widetilde{f})$ is a geodesic, and is clearly independent
of the choice of $\widetilde{f}$, so $b_S (t)$ is indeed well-defined.   
Since developing maps vary smoothly with respect to $t$ we have the following:

\begin{lemma}\label{b:smooth:lemma}
For any $S\in\calS_{\rm gd}\cup\calS_{\rm par}$ the map $t\mapsto b_S (t)$
is smooth with respect to $t$.  
\end{lemma}
 
\subsection{Lifting geometric tetrahedra}
Let now $\Delta$ be a tetrahedron in $\calT$, and denote by
$\widetilde{\Delta}$ a lift of $\Delta$ to $\widetilde{N}'$.
Let $S_i$ be the boundary component of the universal covering
of ${\overline{N}'}$
that corresponds to $v_i$, where $v_1,\ldots,v_4$ are the vertices of
$\widetilde{\Delta}$.
If $\overline{\Delta} (t)$ is the convex hull of $b_{S_1}(t),\ldots,b_{S_4}(t)$
in $\Minkos$, then projecting $\overline{\Delta}(t)$ to
$\matH^3$ and truncating its infinite-volume ends 
with the corresponding $\hat{S}_i (t)$'s gives back a truncated tetrahedron
$\Delta^{\!\ast} (t)$ isometric to the geometric realization 
of $\Delta$ parameterized by $\alpha(t)$.
It is easily seen that the dihedral angles of $\Delta^{\!\ast} (t)$ 
smoothly depend on $b_{S_1}(t),\ldots,b_{S_4}(t)$, so 
\fullref{smoothness:teo} is now reduced to the following:

\begin{prop}\label{smoothness:prop}
If $\dot{\alpha}(0)\in{\rm Ker}\,  dG_{x_0}$, 
then
we can choose the $D_t$'s in such a way that
$\dot{b}_S (0)=0$ for all $S\in\calS_{\rm gd}\cup\calS_{\rm par}$.
\end{prop}

\subsection{The tangent vector to $\rho_t$}
Let us consider the double $DN'$ of $N'$ obtained by mirroring
$N'$ along its boundary (which is now given by some punctured surfaces
of negative Euler characteristic).
Since $\partial N'$ is totally geodesic with respect to the 
hyperbolic structure $M_t$, this structure can be doubled
to a smooth hyperbolic metric $DM_t$ on $DN'$.

Let $D\rho_t\co \pi_1(DN')\to {\rm PSL}(2,\matC)$ be a holonomy representation
corresponding to a developing map for $DM_t$. It is easily seen that we can assume
$D\rho_t$ to vary smoothly with respect to $t$.
We set 
\[
\dot{D\rho}\co \pi_1(DN')\to\mathfrak{sl} (2,\matC),\quad
\dot{D\rho}(\gamma)=\frac{d}{dt}\Big|_{t=0} \left(D\rho_t(\gamma)D\rho_0
(\gamma)^{-1}\right).
\]
As explained in \fullref{cohomology:app},
we have $\dot{D\rho}\in Z^1(\pi_1(DN');\mathfrak{sl} (2,\matC);
 D\rho_0)$. This subsection is devoted to the proof of the following:

\begin{prop}\label{null:cycle:prop}
If $\dot{\alpha}(0)\in{\rm Ker}\,  dG_{x_0}$, then
$\dot{D\rho}\in B^1(\pi_1(DN');\mathfrak{sl} (2,\matC);
 D\rho_0)$. 
\end{prop}

From now on we suppose $\dot{\alpha}(0)\in {\rm Ker}\  dG_{x_0}$.
Let $f_1,\ldots,f_m\in DN'$ be the doubles of the edges 
of $\calT$ and for all $j=1,\ldots,m$
let $\ell_j$ be a small loop in $DN'$ encircling $f_j$. We denote by $\gamma_j$
an element in $\pi_1(DN')$ representing $\ell_j$ (such a $\gamma_j$ is 
well-defined only up to conjugation). Now our hypothesis
implies that if $\theta_j (t)$ is the cone angle of  
(the completion of)
$DM_t$ along $f_j$, then $\dot{\theta}_j(0)=0$. Also observe that we have
\[
D\rho_t(\gamma_j)=g_t\cdot
\left[\left(\begin{array}{cc}
\exp (i\theta_j (t)/2) & 0\\
0 & \exp (-i\theta_j (t)/2)
\end{array}\right)\right]
\cdot g_t^{-1}
\]
for some smooth path $g\co (-\varepsilon,\varepsilon)\to {\rm PSL}(2,\matC)$.
Differentiating this relation we easily get $\dot{D\rho}(\gamma_j)=0$.
Let $K$ be the kernel of the 
map $\pi_1(DN')\to \pi_1(DN)$ induced by the inclusion, and observe that
$K$ is the smallest normal subgroup of $\pi_1 (DN')$
generated by $\gamma_1,\ldots,
\gamma_m$.  Let $D\overline{\rho}_0\co \pi_1 (DN)
\to {\rm PSL}(2,\matC)$ be the natural representation 
associated to $D\rho_0$.
Since
$D\rho_0(\gamma)=1$ and 
$\dot{D\rho}(\gamma)=0$ for all $\gamma\in K$,
\fullref{cohomology1:lem} implies
the following:

\begin{prop}\label{cycle:prop}
A cocycle $z_{{\rho}}\in
Z^1(\pi_1(DN);\mathfrak{sl} (2,\matC);
 D\overline{\rho}_0)$ is naturally induced by
$\dot{D\rho}$.
Moreover, 
$\dot{D\rho}$ belongs to $B^1(\pi_1(DN');\mathfrak{sl} (2,\matC);
 D{\rho}_0)$ 
if and only if $z_\rho$ belongs to 
$B^1(\pi_1(DN);\mathfrak{sl} (2,\matC);
 D\overline{\rho}_0)$.
\end{prop}

Now if $N$ is compact, ie, if no cusps are involved,
\fullref{Wei:teo} directly applies 
concluding the proof of \fullref{null:cycle:prop}.
When there are cusps, 
\fullref{null:cycle:prop}
can be deduced from
\fullref{gar:teo} and the following:

\begin{prop}\label{garland:prop}
If $N$ has cusps, then
$[z_{{\rho}}]\in 
H^1_{\rm par}(\pi_1(DN);\mathfrak{sl}(2,\matC);
D\overline{\rho}_0)$.
\end{prop}
\begin{proof}
Let ${\gamma}\in\pi_1(DN)$ be such that $D\overline{\rho}_0
({\gamma})$ is non-trivial parabolic, and let $\langle
{\gamma}\rangle$ be  
the infinite cyclic
group generated by ${\gamma}$.
We have to check that $z_\rho$ restricts to a coboundary in
$B^1 (\langle {\gamma}\rangle;\mathfrak{sl}(2,\matC);
D\overline{\rho}_0\circ i)$, where $i\co \langle \overline{\gamma}\rangle\to
\pi_1(DN)$ is the natural inclusion.

Without loss of generality
we can suppose
${\gamma}\in \pi_1(T_i)\subset \pi_1(N)\subset \pi_1 (DN)$
for some $i=1,\ldots,k$.
Recall that a preferred
element $\mu_i\in \pi_1 (T_i)$ was previously fixed,
set $u_i (t)=u_i (\alpha (t))$ and observe
that since
$\dot{\alpha}(0)\in {\rm Ker}\  dG_{x_0}$
we have $\dot{u}_i (0)=0$. 
Let ${\mu}'_i,\gamma'$ be elements in $\pi_1 (T_i')$
projecting respectively to $\mu_i,{\gamma}$.
By \fullref{b:smooth:lemma} a smooth
path $g\co (-\varepsilon,\varepsilon)
\to {\rm PSL}(2,\matC)$ exists such that both $g^{-1}_t\cdot
\rho_t({\mu}'_i)\cdot  g_t$ and $g^{-1}_t\cdot
\rho_t(\gamma')\cdot  g_t$
fixes $\infty\in\partial \matH^3$ for  $t\in(-\varepsilon,\varepsilon)$,
so that
\[
\begin{array}{ccc}
D{\rho}_t({\mu}'_i)&=&g_t\cdot
\left[\left(\begin{array}{cc}
\exp (u_i (t)/2) & \tau_i (t)\\
0 & \exp (-u_i (t)/2)
\end{array}\right)\right]\cdot g_t^{-1}\\
D{\rho}_t ({\gamma}')&=&g_t \cdot \left[\left(
\begin{array}{cc}
{a} (t) & {b}(t)\\
0 & {a} (t)^{-1} 
\end{array}\right)\right]
\cdot g_t^{-1}
\end{array}
\]
where $\tau_i,a,b\co (-\varepsilon,\varepsilon)\to\matC$ 
are smooth arcs with $a(0)=1,b(0)\neq 0, \tau_i (0)\neq 0$.

Since $\rho_0 (\pi_1 (T'_i))\cong\matZ+\matZ$
is Abelian, an element $k\in K$ exists such that
${\mu}'_i\cdot \gamma'=k\cdot\gamma'\cdot{\mu}'_i$
in $\pi_1(DN')$. 
By \fullref{cycle:prop},
this readily implies
\[
\frac{d}{dt}\Big|_{t=0} \left(D\rho_t({\mu}'_i) D\rho_t(\gamma')\right)=
\frac{d}{dt}\Big|_{t=0} \left(D\rho_t(\gamma') D\rho_t({\mu}'_i)\right)
\]
which after some computations gives
$\dot{a} (0)=0$ (here we use $\dot{u}_i (0)=0$).
Let us consider the deformation
$\varphi_t\co \langle \gamma' \rangle \to {\rm PSL}(2,\matC)$ defined by
$\varphi_t ((\gamma')^n)=D\rho_t ((\gamma')^n)$ for any
$n\in\matZ$. 

We claim that $\dot{\varphi}=0$ in 
$H^1 (\langle \gamma'\rangle;\mathfrak{sl}(2,\matC);
\varphi_0)$.
This will easily give that $z_\rho$ restricts to a coboundary in
$B^1 (\langle {\gamma}\rangle;\mathfrak{sl}(2,\matC);
D\overline{\rho}_0\circ i)$, whence the conclusion.
Since derivatives of conjugated deformations
differ by a coboundary, we can suppose 
\[
\varphi_t (\gamma')=\left[\left(
\begin{array}{cc}
a (t) & b (t)\\
0 & a (t)^{-1} 
\end{array}\right)\right].
\]
Setting 
\[
v=
\left(
\begin{array}{cc}
\frac{\dot{b} (0)}{2b (0)} & 0\\
0 & -\frac{\dot{b} (0)}{2b (0)} 
\end{array}\right)\in\mathfrak{sl}(2,\matC)
\]
an easy computation shows that $\dot{\varphi}(\gamma')=
v-{\rm Ad}(\varphi_0 (\gamma'))(v)$. This readily implies
that $\dot{\varphi}$ is a coboundary.
\end{proof}

\subsection{The final step}
For $\gamma\in\pi_1(DN')$ 
let ${\rm tr}_\gamma\co (-\varepsilon,
\varepsilon)\to\matC$ be the map defined as follows:
${\rm tr}_\gamma (t)={\rm trace} (D\hat{\rho}_t (\gamma))$,
where $t\mapsto D\hat{\rho}_t(\gamma)\in {\rm SL}(2,\matC)$
is a smooth lift of $t\mapsto D{\rho}_t (\gamma)
\in {\rm PSL}(2,\matC)$ (so ${\rm tr}_\gamma$ is well-defined only up to the
sign). Also recall that $D_t$ and $\rho_t$ are respectively 
a developing map and 
a holonomy representation for the hyperbolic structure $M_t$ on $N'$.
The following result 
can be easily deduced from the proof of the previous proposition.

\begin{lemma}\label{trace:lemma}
Let $\gamma\in\pi_1(DN')$ be such that $D\rho_0 (\gamma)$
is non-trivial parabolic.
Then $\ddot{\rm tr}_\gamma (0)=0$.
\end{lemma}

As a consequence of \fullref{null:cycle:prop} and
\fullref{trace:lemma} we obtain the following:

\begin{prop}\label{zero:derivative:prop}
$D_t$ and $\rho_t$ can be chosen in such a way that
\[
\frac{d}{dt}\Big|_{t=0} \rho_t (\gamma)=0
\]
for all $\gamma\in \pi_1(N')$. Moreover, if $\gamma$ is 
a non-trivial element of a peripheral parabolic subgroup
of $\pi_1 (N')$, then $\ddot{\rm tr}_\gamma (0)=0$.
\end{prop}

Let $\gamma\in\pi_1(N')$ be such that
$\rho_0 (\gamma)\neq 1$
and suppose $r_\gamma\co (-\varepsilon,\varepsilon)\to\partial\matH^3$
is a smooth path such that $r_\gamma (t)$ is a fixed point for
${\rho}_t (\gamma)$ for any $t\in (-\varepsilon,\varepsilon)$.  
We want to study how the derivative of $r(t)$ is related
to the derivative of $\rho_t (\gamma)$. We identify  
$\partial \matH^3$ with $\matC\cup\{\infty\}$ and we set
\[
{\rho}_t (\gamma)=
\left[\left(
\begin{array}{cc}
a(t) & b(t)\\
c(t) & d(t)
\end{array}
\right)\right].
\]
Without loss of
generality we can suppose that 
$\rho_0 (\gamma)$ does not fix $\infty$ and that
$r_\gamma (0)=0\in\matC\subset\partial\matH^3$, 
so $b(0)=0$, $c(0)\neq 0$. An easy computation now shows that
\[
r_\gamma (t)= \left(a(t)-d(t) \pm\sqrt{({\rm tr}_\gamma (t)+2)
({\rm tr}_\gamma (t)-2)}\right)/\left(2c(t)\right).
\]
From this formula we can readily deduce the following lemmas.

\begin{lemma}\label{hyp:fix:lemma}
If $\rho_0 (\gamma)$ is non-trivial loxodromic
and $\dot{\rho} (\gamma)=0$,
then $\dot{r}_\gamma (0)=0$.
On the other hand, let
$\rho_0 (\gamma)$ be non-trivial parabolic.
Also assume that 
$\dot{\rho}(\gamma)=0$
and $\ddot{\rm tr}_\gamma (0)=0$.
Then $\dot{r}_\gamma (0)=0$.
\end{lemma}

We can now conclude the proof of \fullref{smoothness:prop}.
Choose $D_t$ and $\rho_t$ as in the statement of 
\fullref{zero:derivative:prop}.
Let $S$ be a boundary component of the universal covering
of ${\overline{N}'}$
which belongs to $\calS_{\rm par}$, denote by $P_S$ the stabilizer of $S$
in $\pi_1 (N')$ and choose an element $\gamma\in P_S$
with $\rho_0 (\gamma)\neq 1$. By construction the projection
of $b_S (t)$ to $\partial \matH^3$
is fixed by $\rho_t (\gamma)$, so \fullref{hyp:fix:lemma} applies
ensuring $\dot{b}_S (0)=0$.

Suppose now that $S$ belongs to $\calS_{\rm gd}$, 
and let $P_S$ be the stabilizer of $S$
in $\pi_1 (N')$. For $\gamma\in P_S$ with 
$\rho_0 (\gamma)\neq 1$ let 
$p_\gamma (t),q_\gamma (t)$ be the fixed points of $\rho_t (\gamma)$
on $\partial \matH^3$. Note that since $\rho_0 (\gamma)$
is loxodromic we can choose $p_\gamma, q_\gamma$
to vary smoothly with respect to $t$, at least in a small neighbourhood
of $0$. This gives $\dot{p}_\gamma=
\dot{q}_\gamma=0$ by \fullref{hyp:fix:lemma}. 
Now a standard result in Kleinian group theory ensures that the set
$\{ p_\gamma (0),q_\gamma (0):\ \gamma\in P_S,\ \rho_0 (\gamma)\neq 1\}$
is dense in the closure at infinity of $D_0 (S)\subset\matH^3$, 
and this easily implies 
$\dot{b}_S (0)=0$.

\section{Dehn filling}
Once the smoothness of $\Omega (\calT)$ at $x_0$ is established,
one can prove Thurston's hyperbolic Dehn filling Theorem
just by following the strategy described in~\cite{NeuZag}.
At this stage, this argument applies only to those
hyperbolic
manifolds with geodesic boundary which admit a good geodesic 
triangulation (but see \fullref{conj:conj}). 
The following result is taken from~\cite{NeuZag}.

\begin{lemma}\label{neuzag:uv:lemma}
Let $j\in\{1,\ldots,k\}$. Then there exists a complex number
$\tau_j$ with non-zero imaginary part
such that if  
$\{y_n\}_{n\in\matN}\subset \Omega (\calT)$ is a sequence with
$\lim_{n\to\infty} y_n=x_0$ and $u_j (y_n)\neq 0$ for every $n\in\matN$, then
$\lim_{n\to\infty} v_j (y_n)/u_j (y_n)=\tau_j$.
\end{lemma}

\subsection{Thurston's hyperbolic Dehn filling Theorem}
Let $U$ be a sufficiently small neighbourhood of $x_0$
in $\Omega (\calT)$ and let $x\in U$. For $j=1,\ldots,k$, we define
the $j$--\emph{Dehn filling
coefficient} $(p_j (x),q_j (x))\in\matR^2\cup\{\infty\}$
as follows: if $u_j (x)=0$, then $(p_j (x),q_j (x))=\infty$;
otherwise, $p_j (x),q_j (x)$ are the unique real solutions
of the equation
\[
p_j (x) u_j (x)+q_j (x) v_j (x)=2\pi i.
\]
(Existence and uniqueness of such solutions
near $x_0$ can be easily deduced from \fullref{smoothness:teo} 
and \fullref{neuzag:uv:lemma}.)

Let us set
\[
d=(d_1,\ldots,d_k)\co U\to \prod_{i=1}^k S^2,\quad
d_j (x)=(p_j(x),q_j(x))\in S^2=\matR^2\cup\{\infty\}.
\]
As a consequence of \fullref{smoothness:teo}
and \fullref{neuzag:uv:lemma}
we have the following:

\begin{teo}\label{essential:bis:teo}
If $U$ is small enough,
the map $d$ 
defines a
diffeomorphism onto
an open neighbourhood of
$(\infty,\ldots,\infty)$ in $S^2\times\cdots\times S^2$.
\end{teo}

For $x\in\Omega (\calT)$ we denote by $N(x)$
the hyperbolic structure induced on $N$ by $x$, and by
$\widehat{N} (x)$ the metric completion of 
$N(x)$. 
Recall that a preferred basis $\mu_i,\lambda_i$ of $H_1 (T_i;\matZ)$
is fixed for every $i=1,\ldots,k$. 
We also set
\begin{eqnarray*}
I\Omega (\calT)&= &\big\{x\in U\subset\Omega (\calT):
\, {\rm each\ Dehn\ filling\ coefficient\ corresponding\ to}\ x
\\& &\quad {\rm is\ equal\ either\ to}\ \infty
\ {\rm or\ to\ a\ pair\ of\ coprime\ integers}\big\}.
\end{eqnarray*}

\begin{teo}\label{Dehn:fill:eq:teo}
If $U$ is sufficiently small and $x$ belongs
to $I\Omega (\calT)\cap U$, then
$\widehat{N}(x)$ admits a complete finite-volume smooth
hyperbolic structure which is obtained by adding to $N(x)$ a closed
geodesic at any cusp with non-infinite Dehn filling coefficient.
From a topological point of view,
$\widehat{N} (x)$ is obtained by Dehn filling the $i^{\rm th}$ cusp
of $N$ along the slope $p_i (x) \mu_i+q_i (x)\lambda_i$ if
$(p_i (x),q_i (x))\neq\infty$, and by leaving the $i^{\rm th}$ cusp
of $N$ unfilled if $(p_i (x),q_i (x))=\infty$, $i=1,\ldots,k$.
\end{teo}
\begin{proof}
See eg, Thurston~\cite{Thu:bibbia}, Neumann--Zagier~\cite{NeuZag},
Benedetti--Petronio~\cite{BenPet:book}, Frigerio~\cite{Fri:tesi}.
\end{proof}

The following proposition will prove useful in the last subsection.

\begin{prop}\label{fund:prop}
Let $X$ be any smooth manifold and let
$g\co \Omega (\calT)\to X$ be a smooth map. Suppose that
there exists a small neighbourhood $U$ 
of $x_0$ in $\Omega (\calT)$ such that
for all
$x,x'\in U\cap I\Omega (\calT)$ with $d(x)=-d (x')$
we have $g(x)=g (x')$. Then
$dg_{x_0}=0$.
\end{prop}
\begin{proof}
Since $d\co U\to\prod_{i=1}^k S^2$ defines a chart around $x_0$, it is sufficient
to observe that $U\cap I\Omega (\calT)$
accumulates to $x_0$ along any direction in
$T_{x_0} \Omega (\calT)$.
\end{proof}

\subsection{Infinitesimal deformations of compact tetrahedra}
We know from Mostow's rigidity Theorem that compact hyperbolic
$3$--manifolds do not admit deformations. The following results seem to suggest
that in the non-compact case deformations take place mostly near the cusps:
even if it has to be affected by any non-trivial deformation, 
the compact core
offers resistance to changing its shape.
More precisely, we now show that deformations of
$\calT$ near the complete structure affect compact tetrahedra only
at the second order. 

Let $f$ be any \emph{compact} internal edge of $\calT$.
For $x\in\Omega (\calT)$ we denote by $\ell^f (x)$ the
length of $f$ with respect to the metric structure defined by $x$.
By \fullref{trigo:lemma}  the map $\ell^f\co \Omega(\calT)\to\matR$
is smooth.

\begin{prop}\label{length:deform:prop}
We have $d\ell^f_{x_0}=0$.
\end{prop}
\begin{proof}
Let $x,x'\in U\cap I\Omega (\calT)$ be such that $d(x)=-d (x')$.
Then the identity of $N$
extends to a homeomorphism between $\widehat{N} (x)$
and $\widehat{N} (x')$. By Mostow--Prasad's rigidity Theorem,
such a homeomorphism is homotopic to an isometry
$\psi\co \widehat{N} (x)\to \widehat{N} (x')$ via a homotopy
which preserves the geodesic boundary (see eg, \cite{Fri}).
For $y\in U\cap I\Omega (\calT)$
let $f(y)\subset N(y)\subset \widehat{N} (y)$
be the geodesic segment corresponding to $f$.
From the above discussion it follows that
$\psi (f(x))$ is homotopic to $f(x')$ relatively
to $\partial \widehat{N} (x')$. Since both $\psi(f(x))$
and $f (x')$ intersect $\partial \widehat{N} (x')$ 
perpendicularly, this easily implies that
$\psi (f (x))=f (x')$, whence $\ell^f (x)=\ell^f (x')$.
Now the conclusion follows from \fullref{fund:prop}.
\end{proof}

The lengths of the boundary edges of a compact lateral hexagon
smoothly depend on the lengths of its internal edges, and
the dihedral angles of a compact truncated tetrahedron
smoothly depend on the lengths of its internal edges.
Thus \fullref{length:deform:prop}
implies the following results.

\begin{cor}
Fix a boundary edge $f$ lying on a \emph{compact} lateral hexagon
of some tetrahedron of $\calT$, and let $\ell^f\co \Omega (\calT)
\to\matR$ be the function which associates to $x\in\Omega (\calT)$
the length of $f$ in the geometric realization
parameterized by $x$. Then $d\ell^f_{x_0}=0$.
\end{cor}

\begin{cor}\label{compact:deform:cor}
Fix an internal edge $f$  
in a \emph{compact} tetrahedron of $\calT$ and let
$a^f\co \Omega(\calT)$ $\to\matR$ 
be the function that
associates to $x\in\Omega(\calT)$ the dihedral angle 
assigned to $f$ by $x$.
Then
$da^f_{x_0}=0$.
\end{cor}

\subsection{Infinitesimal deformations of the geodesic boundary}
Let ${\rm Teich} (\partial N)$ be
the Teichm\"uller space of hyperbolic structures on $\partial N$,
ie, the space of equivalence classes of hyperbolic metrics
on $\partial N$, where two such metrics are considered equivalent if
they are isometric through a diffeomorphism homotopic to the identity
of $\partial N$. 
For $x\in\Omega (\calT)$ we denote by $B(x)\in {\rm Teich}
(\partial N)$ the equivalence class of the hyperbolic structure
induced by $N(x)$ on $\partial N$.
It is well-known that ${\rm Teich} (\partial N)$
admits a structure of differentiable manifold such that
$B\co \Omega (\calT)\to {\rm Teich} (\partial N)$
is smooth. As a consequence of Mostow--Prasad's rigidity Theorem
and of \fullref{fund:prop}
we get the following:

\begin{prop}\label{diff0:prop}
We have $dB_{x_0}=0$.
\end{prop} 

\appendix

\section{Cohomology theory of representations}
\setobjecttype{App}\label{cohomology:app}

\subsection{The tangent space to a representation}
Let  $G$ be a Lie group
with associated Lie algebra $\mathfrak{g}$
and $\Gamma$ be any group, and denote by
$\calR (\Gamma,G)$ the set of representations
of $\Gamma$ in $G$.
We say that a path 
$\{\rho_t\in\calR (\Gamma,G):\, t\in(-\varepsilon,\varepsilon)\}$
is \emph{smooth} if $\rho_t(\gamma)$ is a smooth function
of $t$ for any $\gamma\in\Gamma$.
If $\{\rho_t\}$ is a smooth path of representations, the tangent
vector to the map $t\mapsto \rho_t(\gamma)$ at $0$
gives an element in $T_{\rho_0 (\gamma)} G$. Identifying
this tangent space with $\mathfrak{g}=T_{1} G$ by \emph{right}
translation we get an element $\dot{\rho} (\gamma)$
of the Lie algebra $\mathfrak{g}$:
\[
\dot{\rho} (\gamma)=\frac{d}{dt}\Big|_{t=0} \left(\rho_t(\gamma) \rho_0 (\gamma)^{-1}\right).
\]
Differentiating the homomorphism relation $\rho_t(\gamma_1\gamma_2)=
\rho_t(\gamma_1)\rho_t(\gamma_2)$ we see that $\dot{\rho}\co \Gamma\to
\mathfrak{g}$ satisfies the so-called \emph{cocycle} relation
\begin{equation*}
\dot{\rho}(\gamma_1\gamma_2)=\dot{\rho}(\gamma_1)+
{\rm Ad}({\rho_0}(\gamma_1)) (\dot{\rho}(\gamma_2)),
\end{equation*}
where ${\rm Ad}\co G\to{\rm GL}(\mathfrak{g})$ is the
usual adjoint representation.  

Consider now a \emph{trivial} deformation of $\rho_0$,
ie, let
$t\mapsto g_t$ be a smooth path
in $G$ starting at the identity and set
$\rho_t(\gamma)=g_t \rho_0(\gamma) g^{-1}_t$ for all $\gamma\in\Gamma$.
Then differentiation shows that
\begin{equation*}
\dot{\rho}(\gamma)=\dot{g}-{\rm Ad}(\rho_0 (\gamma))(\dot{g})
\end{equation*}
for any $\gamma\in\Gamma$, where $\dot{g}\in\mathfrak{g}$ is
the tangent vector to $t\mapsto g_t$ at $t=0$.
We now set:
\begin{eqnarray*}
Z^1 (\Gamma;\mathfrak{g};\rho_0)&=&\{c\co \Gamma\to \mathfrak{g}:
\ c(\gamma_1 \gamma_2)=
c(\gamma_1)+{\rm Ad}(\rho_0(\gamma_1))(c(\gamma_2))\}\\
B^1 (\Gamma;\mathfrak{g};\rho_0)&=&\{b\co \Gamma\to \mathfrak{g}:
\ b(\gamma)=m-{\rm Ad}(\rho_0(\gamma))(m)
\ {\rm for\ some}\ m\in \mathfrak{g}\}\\
H^1(\Gamma;\mathfrak{g};\rho_0)&=&
Z^1 (\Gamma;\mathfrak{g};\rho_0)/B^1 (\Gamma;\mathfrak{g};\rho_0)
\end{eqnarray*}
The above discussion shows that $Z^1(\Gamma;\mathfrak{g};
\rho_0)$ corresponds in some sense 
to the tangent space of $\calR (\Gamma,G)$ at
$\rho_0$.
Under this identification the module
$B^1(\Gamma;\mathfrak{g};\rho_0)$
should represent the tangent space to trivial deformations of $\rho_0$, so 
$H^1(\Gamma;\mathfrak{g};\rho_0)$
should give the tangent space of
$\calR(\Gamma,G)/G$ at $[\rho_0]$ (however, this holds true only
in the setting of algebraic schemes).

Let $c\in Z^1 (\Gamma;\mathfrak{g};\rho)$, and 
suppose that $\Gamma_0$ is a normal subgroup of $\Gamma$
such that
$\rho(\gamma)=1_G$, 
$c(\gamma)=0$ for all $\gamma\in\Gamma_0$.
Let also $\overline{\rho}\co \Gamma/\Gamma_0\to G$
be the representation induced by $\rho$.

\begin{lemma}\label{cohomology1:lem}
The map 
$\overline{c}\co \Gamma/\Gamma_0 \to M$ defined by
$\overline{c} ([\gamma])=c(\gamma)$
is well-defined and gives a cocycle $\overline{c}\in
Z^1 (\Gamma/\Gamma_0;\mathfrak{g};\overline{\rho})$.
Moreover, we have $\overline{c}\in B^1(\Gamma/\Gamma_0;\mathfrak{g};
\overline{\rho})$
if and only if $c\in B^1(\Gamma;\mathfrak{g};\rho)$.
\end{lemma}

\subsection{Classical rigidity results}
Let $N$ be a smooth $3$--manifold without boundary 
and suppose $\rho_0\co \pi_1(N)\to
{\rm PSL}(2,\matC)$ is the holonomy representation for a
complete finite-volume hyperbolic
structure on $N$.
The following result is due to Weil~\cite{Weil1},
and can be considered as a local
version of Mostow's  rigidity Theorem for compact hyperbolic
$3$--manifolds.

\begin{teo}\label{Wei:teo}
Suppose $N$ is compact. Then $H^1(\pi_1(N);\mathfrak{sl}(2,\matC);
\rho_0)=0$.
\end{teo}

Suppose now that $N$ compactifies to
a manifold
$\overline{N}$ with non-empty boundary $\partial \overline{N}=
T_1\sqcup\ldots\sqcup T_k$.
In this case
$\rho_0$ admits non-trivial deformations,
so we cannot expect
$H^1(\pi_1(N);\mathfrak{sl}(2,\matC);
\rho_0)$ to be trivial.
If $K$ is a subgroup
of $\pi_1 (N)$, the natural injection
$i_K\co K\to \pi_1(N)$ induces a map on cohomology
\[
i_K^\ast\co H^1(\pi_1(N);\mathfrak{sl}(2,\matC);
\rho_0)\to H^1(K;\mathfrak{sl}(2,\matC);
\rho_0\circ i_K).
\] 
For $\gamma\in\pi_1 (B)$ we denote by $\langle \gamma\rangle\subset
\pi_1 (N)$ the cyclic
subgroup generated by $\gamma$. If $P=\{\gamma\in\pi_1 (N):\,
\rho_0 (\gamma)\ {\rm is\ parabolic}\}$, we set
\[
H^1_{\rm par}(\pi_1(N);\mathfrak{sl}(2,\matC);
\rho_0)=\bigcap_{\gamma\in P} {\rm Ker}\  i_{\langle \gamma\rangle}^\ast
\subset H^1 (\pi_1(N);\mathfrak{sl}(2,\matC);
\rho_0).
\]
The na\"ive correspondence between
$H^1(\pi_1(N); \mathfrak{sl}(2,\matC);
\rho_0)$ and the set of
conjugacy classes of infinitesimal deformations of $\rho_0$ 
restricts to an identification between $H^1_{\rm par}(\pi_1(N);\mathfrak{sl}(2,\matC);
\rho_0)$ and 
the set of classes of  
infinitesimal deformations
through holonomies for \emph{complete} structures on $N$.
Thus the
following result~\cite{Gar,Rag} can be considered  
an infinitesimal version of
Mostow--Prasad's rigidity Theorem for complete finite-volume
hyperbolic $3$--manifolds:

\begin{teo}\label{gar:teo}
We have $H^1_{\rm par}(\pi_1(N);\mathfrak{sl}(2,\matC);
\rho_0)=0$.
\end{teo}

\bibliographystyle{gtart}
\bibliography{link}

\begin{thebibliography}{}
\providecommand\bibmarginpar{\leavevmode\marginpar}
\def\urlstyle#1{{\tt #1}}

\bibitem{BenPet:book}
\textbf{R Benedetti}, \textbf{C Petronio}, \emph{Lectures on hyperbolic
  geometry}, Universitext, Springer, Berlin (1992) \xox{MR}{1219310}

\bibitem{Choi}
\textbf{Y-E Choi}, \href{http://dx.doi.org/10.1016/j.top.2004.02.002}
  {\emph{Positively oriented ideal triangulations on hyperbolic
  three-manifolds}}, Topology 43 (2004) 1345--1371 \xox{MR}{2081429}

\bibitem{Fri:preprint}
\textbf{R Frigerio}, \emph{Similar fillings and isolation of cusps of
  hyperbolic $3$-manifolds}, Technical report \xox{arXiv}{math.GT/0504147}

\bibitem{Fri}
\textbf{R Frigerio}, \emph{Hyperbolic manifolds with geodesic boundary which
  are determined by their fundamental group}, Topology Appl. 145 (2004) 69--81
  \xox{MR}{2100865}

\bibitem{Fri:tesi}
\textbf{R Frigerio}, \emph{Deforming triangulations of hyperbolic $3$-manifolds
  with geodesic boundary}, PhD thesis, Scuola Normale Superiore, Pisa (2005)

\bibitem{FriMaPe2}
\textbf{R Frigerio}, \textbf{B Martelli}, \textbf{C Petronio},
  \href{http://projecteuclid.org/getRecord?id=euclid.em/1090350932}
  {\emph{Small hyperbolic 3-manifolds with geodesic boundary}}, Experiment.
  Math. 13 (2004) 171--184 \xox{MR}{2068891}

\bibitem{FriPe}
\textbf{R Frigerio}, \textbf{C Petronio},
  \href{http://dx.doi.org/10.1090/S0002-9947-03-03378-6} {\emph{Construction
  and recognition of hyperbolic 3-manifolds with geodesic boundary}}, Trans.
  Amer. Math. Soc. 356 (2004) 3243--3282 \xox{MR}{2052949}

\bibitem{Gar}
\textbf{H Garland}, \emph{A rigidity theorem for discrete subgroups}, Trans.
  Amer. Math. Soc. 129 (1967) 1--25 \xox{MR}{0214102}

\bibitem{HodKer}
\textbf{C\,D Hodgson}, \textbf{S\,P Kerckhoff}, \emph{Rigidity of hyperbolic
  cone-manifolds and hyperbolic {D}ehn surgery}, J. Differential Geom. 48
  (1998) 1--59 \xox{MR}{1622600}

\bibitem{Koj:proc}
\textbf{S Kojima}, \emph{Polyhedral decomposition of hyperbolic manifolds with
  boundary}, Proc. Work. Pure Math. 10 (1990) 37--57

\bibitem{Luo}
\textbf{F Luo}, \href{http://dx.doi.org/10.1090/S1079-6762-05-00142-3} {\emph{A
  combinatorial curvature flow for compact 3-manifolds with boundary}},
  Electron. Res. Announc. Amer. Math. Soc. 11 (2005) 12--20 \xox{MR}{2122445}

\bibitem{NeuZag}
\textbf{W\,D Neumann}, \textbf{D Zagier},
  \href{http://dx.doi.org/10.1016/0040-9383(85)90004-7} {\emph{Volumes of
  hyperbolic three-manifolds}}, Topology 24 (1985) 307--332 \xox{MR}{815482}

\bibitem{PetPor}
\textbf{C Petronio}, \textbf{J Porti}, \emph{Negatively oriented ideal
  triangulations and a proof of {T}hurston's hyperbolic {D}ehn filling
  theorem}, Expo. Math. 18 (2000) 1--35 \xox{MR}{1751141}

\bibitem{Rag}
\textbf{M\,S Raghunathan}, \emph{Discrete subgroups of {L}ie groups},
  Ergebnisse series 68, Springer, New York (1972) \xox{MR}{0507234}

\bibitem{Thu:bibbia}
\textbf{W\,P Thurston}, \emph{The geometry and topology of $3$-manifolds}
  (1979)Mimeographed notes

\bibitem{Weil1}
\textbf{A Weil},
  \href{http://links.jstor.org/sici?sici=0003-486X(196407)2:80:1%3C149:ROTCOG%%
3E2.0.CO%3B2-1} {\emph{Remarks on the cohomology of groups}}, Ann. of Math.
  $(2)$ 80 (1964) 149--157 \xox{MR}{0169956}

\end{thebibliography}

\end{document}